\documentclass[12pt]{article}
\usepackage{amsfonts}
\usepackage{amsmath}
\usepackage[english]{babel}
\newtheorem{prop}{Proposition}[section]

\newtheorem{theo}[prop]{Theorem}
\newtheorem{lem}[prop]{Lemma}

\newtheorem{rem}[prop]{Remark}

\def\1{1\!{\rm l}}

\begin{document}
\title{Results dealing with the behavior of the integrated density of states of random divergence operators}
\maketitle
\begin{center}\bf{{Hatem Najar}} \footnote{D\'epartement de
Math\'ematiques Physiques I.P.E.I. Monastir, 5000 Monastir Tunisie
. \\Researches partially supported by  CMCU N 04/S1404 and
Research Unity 01/UR/ 15-01 projects. }
\end{center}
\begin{abstract}

\vskip.5cm\noindent In this paper we generalize and improve results proven for
acoustic operators in \cite{jmp,long}. It deals with the behavior of
the integrated density of states of random divergence operators of the form
$H_\omega=\sum_{i,j=1}^d\partial_{x_{i}}a_{i,j}(\omega,x)\partial_{x_j}$;
 on the internal band edges of the spectrum. We propose an application of
such a result to get localization.

\end{abstract}
{\small\sf 2000 Mathematics Subject Classification :81Q10, 35P05,
37A30,
47F05.}
\setcounter{section}{0}
\section{Introduction}
Let us consider the random divergence operator

\begin{equation}
 \label{b}
 H_\omega=-\nabla
 {\rho_\omega}\nabla=\sum_{i,j=1}^d\partial_{x_{i}}a_{i,j}(\omega,x)\partial_{x_j};
\end{equation}
 where $\rho_\omega=(a_{i,j}(\omega,x))_{1\leq i,j \leq d}$ is an
elliptic, $d\times d$-matrix valued, $\mathbb{Z}^d$-ergodic random
field. i.e there exists some constant $\rho_*>1$, satisfying
\begin{equation}
 \label{eq:3}
 \frac{1}{\rho_*}|\xi|^2\leq \langle \rho_\omega \xi,\xi \rangle\leq\rho_*|\xi|^2,\ \ \forall \xi\in\mathbb{C}^d.
\end{equation}
 This operator describes a vibrating membrane in the random medium as well as in the particular case when
$\rho_{\omega}=\frac{1}{\varrho_{\omega}}\cdot I_{d}$, ($I_d$ is
the identity matrix and $\varrho_{\omega}$ is a real function ) we
get the acoustic operator \cite{Figkl1,jmp,long}. The great
interest of this operator both from the physical and the
mathematical point of view is quite obvious and known \cite{Stob}.
\newline   As this paper is devoted to the study of the behavior
of the integrated density of states, we start by recalling that it
is defined as follows: We note by $H_{\omega,\Lambda }$ the
restriction of $H_{\omega }$ to $\Lambda$ with self-adjoint
boundary conditions. As $H_{\omega }$ is elliptic, the resolvent
of $H_{\omega,\Lambda }$ is compact and, consequently, the
spectrum of $H_{\omega,\Lambda }$ is discrete and is made of
isolated eigenvalues of finite multiplicity \cite{21}. We define
\begin{equation}
N_{\Lambda }(E)=\frac{1}{\mathrm{{vol}(\Lambda )}}\cdot \#\{\mathrm{{%
eigenvalues\ of\ }H_{\omega,\Lambda }\leq E\}.} \label{aha}
\end{equation}
Here $\mathrm{{vol}(\Lambda )}$ is the volume of $\Lambda $ in the
Lebesgue sense and $\#E$ is the cardinal {of $E$.}\newline
It is shown that the limit of $N_{\Lambda }(E)$ when $\Lambda $ tends to $%
\mathbb{R}^{d}$ exists almost surely and is independent of the
boundary
conditions. It is called the \textbf{integrated density of states} of $%
H_{\omega }$ (IDS as an acronym). See \cite{kima7}.\subsection{The
result} The essential goal of this work is to study internal
Lifshitz tails for the operator defined by (\ref{b}). We review
some results proven previously and improve them. In
\cite{jmp,long} we have studied the same question
 under a special  regime of disorder, precisely under the assumption that
$ \displaystyle \lim_{\varepsilon\to 0}\frac{\log|\log
\mathbb{P}\{ \omega_{0}\leq \varepsilon \}|} {\log \varepsilon
}=0$. The main novelty here is that in the present work we omit
this condition and we take a more general distribution of the
random variables precisely we consider the case when
$\displaystyle \log\mathbb{P}\{\omega_{0}\leq \varepsilon\}\sim
-\varepsilon^{-\kappa}$ and we extend the result for another class
of random
 Schr\"odinger operators. \newline
As a possible application of our result we get localization near the
band edges which is based on the fact that near those edges the
integrated density of states exhibits Lifshitz tails.
The main interest of such a technique lies in a much weaker assumption on the probability
distribution.
\subsection{The model}
 Assume that $\rho_\omega$ is of Anderson type i.e.
it has the form
\begin{equation}
 \label{eq:6}
 \rho_\omega(x)=\rho^+(x)+\sum_{\gamma\in\mathbb{Z}^d}\omega_\gamma
 \rho^0(x-\gamma),
\end{equation}where\newline
$\bf{(A.0)}$\\
\begin{itemize}
\item $\rho^+=(\rho_{i,j}^+)_{1\leq i,j\leq d}$ is a, $\mathbb{Z}^d$-periodic and elliptic $d\times d$-matrix
valued function.
\item $\rho^0=(\rho_{i,j}^{0})_{1\leq i,j\leq d}$ is a $d\times
d$-matrix valued function such that for some elliptic matrix,
$\varrho_{+}^{0}$, we have
$$0\leq \sum_{\gamma\in \mathbb{Z}^d}\rho^{0}(x-\gamma)\leq \varrho_{+}^{0}.$$
\item $(\omega_\gamma)_{\gamma\in\mathbb{Z}^d}$ is a family of non constant
and positive, independent identically distributed random variables
taking values in $[0,1]$.We note by
$(\Omega,\mathcal{F},\mathbb{P})$ the probability space and we
suppose that
 \begin{equation}
 \lim_{\varepsilon \to 0^+}\frac{\log|\log\mathbb{P}(\{\omega_{0}\leq \varepsilon\})|}{\log \varepsilon}=-\kappa,\ \ \ \kappa\in [0,+\infty[.
 \label{you1}\end{equation}
 \end{itemize}
 Let $h(\rho_{\omega})$ be the quadratic form defined as follows: For $%
u\in H^{1}({\mathbb{R}}^{d})=\mathcal{D}({\textit{h}}(\rho
_{\omega }))$
\begin{eqnarray*}
h(\rho_{\omega })[u,u] &=&\int_{{\mathbb{R}}^{d}}\rho
_{\omega}
(x) \nabla u(x)\cdot\overline{\nabla u(x)}dx \\
&=&\sum_{1\leq i,j\leq d}\int_{{\mathbb{R}}^{d}}a_{i,j}(\omega,x)
\partial _{x_{i}}u(x)\overline{\partial _{x_{j}}u(x)}dx.
\end{eqnarray*}
$h(\rho_{\omega })$ is a positive and closed quadratic form.
$H_{\omega }$ given by (\ref{b}) is defined as the self adjoint
operator associated to ${\textit{h}}(\rho_{\omega })$ \cite{21}.
By this, $H_{\omega}$ is a measurable family of self adjoint and
ergodic operators\cite{Figkl1,Stob}.
\subsection{Reference operator}
It is convenient to write $H_{\omega }$ as a perturbation of some
background periodic operator $H_{0}$. More precisely we write:
\begin{equation*}
H_{\omega }=H_{0}+V_{\omega },
\end{equation*}
with
\begin{equation*}
H_{0}=-\nabla {\rho^+}\nabla \end{equation*} and
\begin{equation*}V_{\omega }(\cdot)=-\nabla (\sum_{\gamma\in
\mathbb{Z}^d}\omega_{\gamma}{\rho^0}(\cdot-\gamma))\nabla\geq 0.
\end{equation*}
\subsubsection{Some facts from Floquet theory} Now we review some standard
facts from the Floquet theory for periodic operators. Basic
references of this material can be found in \cite{21}.\newline As
$\rho^+$ is a ${{\mathbb{Z}^{d}}}$-periodic
matrix, for any ${\gamma \in {\mathbb{Z}%
^{d}}}$, we have
\begin{equation*}
\tau _{\gamma }H_{0}\tau _{\gamma }^{\ast }=\tau _{\gamma
}H_{0}\tau _{-\gamma }=H_{0}.
\end{equation*}
Let $\mathbb{T}^{\ast }=\mathbb{R}^{d}/(2\pi {\mathbb{Z}^{d}})$. We define ${%
{\mathcal{H}}}$ by
\begin{equation*}
{{\mathcal{H}}}=\{u(x,\theta )\in L_{loc}^{2}({\mathbb{R}}^{d})\otimes L^{2}(%
{\mathbb{T}}^{\ast });\forall (x,\theta,\gamma )\in
\mathbb{R}^{d}\times \mathbb{T}^{\ast }\times \mathbb{Z}^{d};\
u(x+\gamma,\theta )=e^{i\gamma \theta }u(x,\theta )\}.
\end{equation*}
There exists $U$ a unitary isometry from $L^{2}({\mathbb{R}}^{d})$ to ${{%
\mathcal{H}}}$ such that $H_{0}$ admits the Floquet decomposition
\cite{21}
\begin{equation*}
UH_{0}U^{\ast }=\int_{{\mathbb{T}^{\ast }}}^{\oplus }H_{0}(\theta
)d\theta.
\end{equation*}
Here $H_{0}(\theta )$ is the operator $H_{0}$ acting on ${{\mathcal{H}}}%
_{\theta }$, defined by
\begin{equation*}
{\mathcal{H}}_{\theta }=\{u\in L_{loc}^{2}(\mathbb{R}^{d});\forall
\gamma \in {\mathbb{Z}^{d}},u(x+\gamma )=e^{i\gamma \theta
}u(x)\}.
\end{equation*}
As $H_{0}$ is elliptic, we know that, $H_{0}(\theta )$ has a
compact resolvent; hence its spectrum is discrete \cite{21}. We
denote its eigenvalues, called Floquet eigenvalues of $H_{0}$, by
\begin{equation*}
E_{0}(\theta )\leq E_{1}(\theta )\leq \cdot \cdot \cdot \leq
E_{n}(\theta )\leq \cdot \cdot \cdot.
\end{equation*}
The functions $(\theta \rightarrow E_{n}(\theta ))_{n\in {%
\mathbb{N}}}$ are Lipshitz-continuous, and we have
\begin{equation*}
E_{n}(\theta )\rightarrow +\infty \ \ \mathrm{as}\ n\rightarrow
+\infty \ \ \mathrm{uniformly\ in}\ \ \theta.
\end{equation*}
The spectrum $\sigma (H_{0})$ of $H_{0}$ has a band structure.
(i.e $\displaystyle \sigma (H_{0})=\cup _{n\in
\mathbb{N}}E_{n}({\mathbb{T}^{\ast }}).$)
\newline
The periodic operator $H_{0}$ has an IDS which will be denoted by
$n$. The behavior of $n$ at a band edge $E_{+}$, is said to be
{\bf{non-degenerate}} if,
\begin{equation}
\lim_{\varepsilon\to 0^+}\frac{\log |n(E_{+}+\varepsilon)-n(E_{+})|}{%
\log\varepsilon}=\frac{d}{2}. \label{tun}
\end{equation}
\subsubsection{The main
assumptions} As we study internal Lifshitz tails it is naturel to assume that $H_{0}$ has a spectral gap below
$E_{+}$. More precisely we assume that:
\newline $\mathbf{(A.1)}$ \newline There exists $E_{+}$ and
$\delta >0$ such that $\sigma (H_{0})\cap \lbrack
E_{+},E_{+}+\delta )=[E_{+},E_{+}+\delta )$ and $\sigma
(H_{0})\cap (E_{+}-\delta,E_{+}]=\emptyset $. \newline
As, $V_{\omega }\geq 0$, the spectrum $\Sigma $ of $H_{\omega}$ contains an
interval of the form $%
[E_{+},E_{+}+a]\ (a>0)$ \cite{Kirmar}. As we are interested in the behavior of
the IDS in the neighborhood of $E_{+} $, we require that $E_{+}$
remains always the edge of a gap for $\Sigma $,
when the perturbation is turned on. More precisely, if for all $%
t\in \lbrack 0,1]$, we define $H_{\omega,t}=H_{0}+tV_{\omega }$ and $%
\Sigma _{t}$ is the almost sure spectrum of $H_{\omega,t}$, then
one requires that the following assumption holds. \newline
$\mathbf{(A.2)}$ \newline There exists $\delta ^{\prime }>0$ such
that for all $t\in \lbrack 0,1],\ \Sigma _{t}\cap \lbrack
E_{+}-\delta ^{\prime },E_{+})=\emptyset $.\newline We assume also
the following:\newline $\bf{(A.3)}$ \\ We will state that the
behavior of the IDS depends on the form of the
perturbation. One distinguishes between two behaviors of $\rho^{0}$. \newline Let $C_{0}= \{x\in \Bbb{R}%
^d; \forall 1\leq j \leq d; -\frac{1}{2}<x_j\leq \frac{1}{2}\}$
and let $0< g_{-}<g_{+}$ be two positives constants. \newline 1.
$\rho^0$ is of long range type.
\\ There exists $\nu \in (d,d+2]$ such that for any $\gamma \in
\mathbb{Z}^{d}$, $1\leq i,j\leq d$ and almost every $x$ in $C_{0}$
one has
\begin{equation}
g_{-}\leq \rho_{i,j}^0(x-\gamma)\cdot (1+|\gamma |)^{\nu }\leq
g_{+}; \label{f}
\end{equation}
and
\begin{equation}
g_{-}\leq |\partial{x_{i}}(\rho_{i,j}^0)(x-\gamma)|\cdot
(1+|\gamma |)^{\nu }\leq g_{+}.\label{f2}
\end{equation}
2. $\rho^0$ is of short range type.\newline There exists $\nu
>d+2$ such that for any $\gamma \in \mathbb{Z}^{d}$, $1\leq
i,j\leq d$ and almost every $x$ in $C_{0}$ one has
\begin{equation}
0\leq \rho_{i,j}^0(x-\gamma)\cdot (1+|\gamma |)^{\nu }\leq
g_{+}; \label{f}
\end{equation}
and
\begin{equation}
0\leq |\partial{x_{i}}(\rho_{i,j}^0)(x-\gamma)|\cdot (1+|\gamma
|)^{\nu }\leq g_{+}.\label{f2}
\end{equation}
\section{Results and discussions }
The main result of this paper is the following:
\begin{theo}\label{T0}
Let $H_{\omega}$ be the operator defined by (\ref{b}). We assume
that {\bf{(A.1)}}, {\bf{(A.2)}} hold. Then if\newline 1. $\rho^0$
is of long range type then,
\begin{multline*}
\lim_{\varepsilon \to 0^{+}} \frac{\log\big(n(E_++\varepsilon)-n(E_+)\big)}{\log\varepsilon} =\frac{d}{2}\Rightarrow
\\ \lim_{\varepsilon \to 0^+} \frac{\log|\log \big( N(E_+
+\varepsilon)-N(E_+)\big)|}{\log\varepsilon}=-\sup(\frac{d}{2}+\kappa,\frac{d}{\nu-d}),
\end{multline*}if $\kappa +\frac{d}{2}<\frac{d}{\nu-d}$,
\begin{equation}\lim_{\varepsilon \to 0^+} \frac{\log|\log \big( N(E_+
+\varepsilon)-N(E_+)\big)|}{\log\varepsilon}=-\frac{d}{\nu-d}.
\end{equation}
2. $\rho^0$ is of short range type then,
\begin{multline*}
\lim_{\varepsilon \to 0^+} \frac{\log|\log \big( N(E_+
+\varepsilon)-N(E_+)\big)|}{\log\varepsilon}=-(\frac{d}{2}+\kappa)\Leftrightarrow \\
\lim_{\varepsilon \to 0^{+}} \frac{\log\big(n(E_++\varepsilon)-n(E_+)\big)}{%
\log\varepsilon} =\frac{d}{2}.
\end{multline*}
\end{theo}
\begin{rem}
\label{R00} The result of Theorem \ref{T0} is stated for lower
band edges. Under adequate assumptions the corresponding result
is true for upper band edges.
\end{rem}
Now, let us comment the result. According to the Theorem \ref{T0}
one notices that the behavior of the random variables is linked up
to the lifshitz exponent, and determines if one is located in a
classical regime or in a quantum one; i.e   if the kinetic energy
intervenes or if it does not in the Lifshitz exponent. In the long
range case one sees that it depends on the value of $\kappa$, the
Lifshitz asymptotics are classical (if
$\kappa<\frac{d}{\nu-d}-\frac{d}{2}$) or quantum (if
$\kappa>\frac{d}{\nu-d}-\frac{d}{2}$). In other terms in the case
of the long range potential, Lifshitz exponent depends on the
uncertainty principle, i.e on the kinetic energy only in the case
when ($\frac{d}{\nu -d}<\kappa+\frac{d}{2}$). In contrast, when
($\frac{d}{\nu-d}>\kappa+\frac{d}{2}$) then the Lifshitz
asymptotics are not governed by the same considerations. This is
due to the fact that in the long range case as the potential
decreases slowly, locally the potential is an empirical average of
random variables. This leads to the fact that its effect is more
important and more influencing than the spatial extension of the
considered state.\newline From what it has been said previously,
one concludes that the value of $\kappa$ is responsible for the
transition between those two regimes.\newline The proof of the
main result is now classic and based on the technic of periodic
approximations which where originally stated by Klopp in
\cite{kp8}. It is quite close and follows the same steps used in
\cite{jmp,long}. We omit details and we refer the reader to the
above references.
\subsection{Application}
Now we state a useful result which can be related to the Theorem \ref{T0}. Let
\begin{theo}\label{2d}
 Let $\theta \in \mathbb{R}^d$ and $E_{+}>0$ a band edge of the
 spectrum of $H_{\omega}$. Then for any $\alpha >1$, integer $p>0$, for
$k\in \mathbb{N}$ sufficiently large,
 one has
 $${\bf{{(P1)}}}\
\mathbb{P}\big(\big\{dist\big(\sigma(H_{\omega,\Lambda_{k^{\alpha}}}^{\theta}),E_{+}\big)\leq
 \frac{1}{k}\big\}\big)\leq \frac{1}{k^{p}}.
 $$
  Where $\Lambda_{k}$ is the box centered in $0$ of side length
 $2k+1$ and $A_{\omega,\Lambda_{k}}^{\theta}$ is the operator
 $H_{\omega}$ restricted to this box with $\theta$-quasiperiodic
 boundary condition i.e with boundary condition
 $\varphi(x+\gamma)=e^{i\gamma\cdot \theta}\varphi (x)$ for any
 $\gamma\in (2 k+1)\mathbb{Z}^d$.
 \end{theo}
To be able to apply the multiscale analysis \cite{Figkl1,Stob}, we
assume that $\rho^{0}$ is compactly supported. Indeed, when the single
site is compactly supported $H_{\omega}$ satisfies a Wegner
estimate \cite{Figkl1} i.e for some $\alpha
>0$ and $n>0$ for $E\in \mathbb{R}$ for $k\geq 1$ and
$0<\varepsilon <1 $, there exists $C(E)>0$ such that one has
\begin{equation}
{\bf{{(P2)}}}\
\mathbb{P}\big(\big\{dist(\sigma(H_{\omega,\Lambda_{k}}^{\theta}),E)\leq
 \varepsilon \big\}\big)\leq C(E)\cdot |\Lambda_{k}|^{\alpha}\cdot
 \varepsilon^{n}.\label{w}
\end{equation}
So, for a band edge energy $E_{+}$ using the Theorem \ref{2d} for
$\theta =0$, we obtain the initial estimate to start a multi-scale
analysis. This proves that the spectrum of $H_{\omega}$ is
exponentially localized in some interval around the energy $E_{+}$
i.e that in some neighborhood of $E_{+}$ eigenfunctions associated
to energies in that interval are exponentially localized. More
precisely we have
\begin{theo}
\label{CAF1} Let $H_{\omega}$ defined by (\ref{b}).
 We assume that
 {\bf{(A.1)}} and {\bf{(A.2)}} hold and the single site is compactly
supported. There exists $\varepsilon_{0}>0$ such that\\ (i)
$\Sigma \cap [E_+,E_++\varepsilon_{0}]= \Sigma _{pp}\cap
[E_{+},E_{+}+\varepsilon_{0}]$.\\
(ii) an eigenfunction corresponding to an eigenvalue in
{${[E_+,E_++\varepsilon_{0}]}$} decays exponentially.\\ (iii) for
all  $p>0$,
$$\mathbb{E}\Big\{\sup_{t>0}\Big|\Big|\Big|X\Big|^p
e^{itH_{\omega}}
P_{[E_+,E_++\varepsilon_{0}]}(H_{\omega})\chi_{K}\Big|\Big|\Big\}<+\infty.
$$ Here $P_{I}(H_{\omega})$ is the spectral projection on
 the interval $I,$ $\chi_{K}$ is the characteristic function of
$K$, $K$ is a compact of $\mathbb{R}^d $ and $X$ is the position
operator.
\end{theo}
To comment upon Theorem \ref{CAF1}, let us consider the wave
equation:
\begin{equation}\frac{\partial ^{2}u}{\partial
t^2}=H_{\omega}u.\label{sa}
\end{equation} The solution of (\ref{sa}) is given \cite{Stob} by
$$
u(t,\cdot)=\cos(t\sqrt{H_{\omega}})u_0+\sin(t\sqrt{H_{\omega}})u_1,
$$

where $u_0=u(0,\cdot)$ and
$\sqrt{H_{\omega}}u_1=(\partial_{t}u)(0,\cdot)$ denote the initial
data.\newline The result of Theorem \ref{2d} and the one of
Theorem \ref{CAF1} can be related to the behavior of the
integrated density of states in the neighborhood of the so-called
fluctuation boundary $E_{+}$\cite{kp5,Pasfig}. This is done in the
Schr\"odinger case in \cite{Ves}.
\subsection{The periodic approximations}\label{abour}
Let us consider the following periodic operator
\begin{equation*}
H_{\omega,k}=-\nabla \cdot \rho_{\omega,k}\cdot \nabla,
\end{equation*}
where $\rho_{\omega,k}$ is the following matrix
\begin{equation*}
\rho_{\omega,k}=\rho^++\sum_{\gamma \in C_{k}\cap
{\mathbb{Z}^{d}} }\omega _{\gamma }\sum_{\beta \in
(2k+1){\mathbb{Z}^{d}}}\rho^{0}(\cdot -\gamma -\beta).
\end{equation*}
$C_{k}$ is the cube
\begin{equation*}
C_{k}=\{x\in {\mathbb{R}}^{d};\forall 1\leq j\leq d,\ -\frac{2k+1}{2}%
<x_{j}\leq \frac{2k+1}{2}\}.
\end{equation*}
We set
$$
V_{\omega,k}=-\nabla \cdot \Big( \sum_{\gamma \in C_{k}\cap
{\mathbb{Z}^{d}} }\omega _{\gamma }\sum_{\beta \in
(2k+1){\mathbb{Z}^{d}}}\rho^{0}(\cdot -\gamma -\beta)\Big)\cdot \nabla
$$
\newline
$H_{\omega,k}$ is $(2k+1){\mathbb{Z}^{d}}$-periodic and
essentially self
adjoint operator. Let ${\mathbb{T}}_{k}^{\ast }=({\mathbb{R}}%
^{d})/2(2k+1)\pi {\mathbb{Z}^{d}}$. We define $N_{\omega,k}$ the IDS of $%
H_{\omega,k}$ by
\begin{equation}
N_{\omega,k}(E)=\frac{1}{(2\pi )^{d}}\sum_{n\in
\mathbb{N}}\int_{\{\theta \in {\mathbb{T}}_{k}^{\ast },\ E_{\omega
,k,n}(\theta )\leq E\}}d\theta. \label{11}
\end{equation}
Let $dN_{\omega,k}$ be the derivative of $N_{\omega,k}$ in the
distribution sense. As $N_{\omega,k}$ is increasing, $dN_{\omega
,k}$ is a positive measure; it is the density of states of
$H_{\omega,k}$. We denote by $dN$ the density of states of
$H_{\omega }$.
\begin{theo} \label{T20}\cite{kp1,jmp}
For any $\varphi\in \Lambda_{0}^{\infty}({\Bbb{R}})$ and for almost all
$\omega\in \Omega $ we have
$$\lim_{k\rightarrow \infty
}\langle \varphi,dN_{\omega,k}\rangle=\langle \varphi,dN\rangle.$$
\end{theo}
In what follows we give a well-known \cite{kp8,jmp, long} result stating that the IDS of
$H_{\omega}$ is exponentially well approximated by the expectation
of the IDS of the periodic operators $H_{\omega,k}$ when $k$ is
polynomial in $\varepsilon^{-1}$. More precisely
\begin{lem}\label{mez}
For any $\eta_{0}>1$ and $I\subset\mathbb{R}$ a compact interval,
there exists $\nu_{0}>0$ and $\varepsilon_{0}>0$ such that, for
$0<\varepsilon<\varepsilon_{0},E\in I$ and $k\geq k_{1}=
\varepsilon^{-\nu_{0}}$, we have
\begin{multline}
\label{sammes}
\mathbb{E}[N_{\omega,k}(E+\varepsilon/2)-N_{\omega,k}(E-\varepsilon/2)]-e^{\varepsilon^{-\eta_{0}}}\\
\leq N(E+\varepsilon)-N(E) \\ \leq
\mathbb{E}[N_{\omega,k}(E+2\varepsilon)-N_{\omega,k}
(E-2\varepsilon)]+e^{-\varepsilon^{-\eta_{0}}}.
\end{multline}
\end{lem}
\section{The proof of Theorem \ref{T0}}
To prove Theorem \ref{T0}, we use periodic
approximations. We prove a lower and an upper bounds on $N(E_{+}+%
\varepsilon)-N(E_{+})$. The upper and lower bounds are proven
separately.
\subsection{The lower bound}
We postpone the proof of the lower bound. More details can be
found in \cite{jmp,long}. It consists in proving the following
theorem.
\begin{theo}\label{LL1} Let $H_{\omega}$ be the operator defined by (\ref{b}).
 We assume that
{\bf{(A.1)}}, {\bf{(A.2)}} hold. Then, \\
$\bullet$ if $\rho^{0}$ is of long range type, we have
\begin{equation} \liminf_{\varepsilon \to 0^{+}}
\frac{\log\Big|\log \Big ( N(E_+
 +\varepsilon)-N(E_+)\Big)\Big|}{\log\varepsilon}\geq-
\sup(\frac{d}{2}+\kappa,\frac{d}{\nu-2}).\label{hifr}
\end{equation}
$\bullet$ if $\rho^{0}$ is of short range type, we have
\begin{equation} \liminf_{\varepsilon \to 0^{+}}
\frac{\log\Big|\log \Big ( N(E_+
 +\varepsilon)-N(E_+)\Big)\Big|}{\log\varepsilon}\geq
- (s\frac{d}{2}+\kappa).\label{mouch1}
\end{equation}Here $s<1$ if $n$ is degenerate and $s=1$ if not.\end{theo}
\textbf{Proof.} By assumption, there is a spectral gap below $E_+$
of length at least $\delta ^{\prime}> 0$. Thus, for
$\varepsilon<\delta ^{\prime}$ we have
$$
N(E_++\varepsilon) - N(E_+)=N(E_++\varepsilon) -
N(E_+-\varepsilon).
$$
To prove Theorem \ref{LL1}, it suffice to lower bound $N(E_++\varepsilon)-N(E_+-%
\varepsilon)$. Then, for $N$ large, we will show that $H_{\omega,
\Lambda_N}$ ($H_{\omega, \Lambda_{N}}$ is $H_{\omega}$ restricted
to $\Lambda_{N}$ with
Dirichlet boundary conditions) has a large number of eigenvalues in $%
[E_+-\varepsilon,E_++\varepsilon]$ with a large probability. For
this we will construct a family of approximate eigenvectors
associated to approximate eigenvalues of $H_{\omega, \Lambda_N}$
in $[E_+-\varepsilon, E_++\varepsilon]$. These functions can be
constructed from an eigenvector of $H_{0}$ associated with $E_+$.
Locating this eigenvector in $\theta $ and imposing to
$\omega_{\gamma}$ to be small for $\gamma$ in some well chosen
cube, one obtains an approximate eigenfunction of $H_{\omega,
\Lambda_N}$. Locating the eigenfunction in $x$ in several
disjointed places, we get several eigenfunctions two by two
orthogonal.The subtlety  is in the good choice of the size of the cube.
\newline
Using the same computation done in \cite{jmp,long} we get that we
have to estimate the following two probabilities: For
$1>\alpha>0$,
$$\mathbb{P}_{\varepsilon,\alpha,1}=\mathbb{P}\Big(\Big\{\omega;\ |\beta|\leq \varepsilon^{-(1+\alpha)/2};\ \sum_{\gamma\in \mathbb{Z}^d}
\omega_{\gamma}(1+|\beta-\gamma|)^{-\nu}\leq \varepsilon ^{1+\alpha}\Big\}\Big),$$
$$
\mathbb{P}_{\varepsilon,\alpha,2}=\mathbb{P}\Big(\Big\{\omega;\ \sum_{\gamma\in \Lambda_{\alpha}
(\varepsilon^{s})}\omega_{\gamma}(1+|\gamma|)^{-\nu}\leq \frac{\varepsilon^{1+\alpha}}{2}\Big\}\Big);
$$
here $\Lambda_{\alpha}(\zeta)=\{\gamma\in \mathbb{Z}^d;\ \forall
1\leq j\leq d; \ |\gamma_j|\leq
\zeta^{-(\frac{1}{2}+\alpha)}\}$.\newline Indeed we have the
following relations\cite{jmp,long}:\newline $\bullet$ if
$\rho^{0}$ is of long range type, we have
\begin{equation} \liminf_{\varepsilon \to 0^{+}}
\frac{\log\Big|\log \Big
( N(E_++\varepsilon)-N(E_+)\Big)\Big|}{\log\varepsilon}\geq \liminf_{\varepsilon \to 0^{+}}
 \frac{\log\Big|\log ( \mathbb{P}_{\varepsilon,\alpha,1})\Big|}{\log\varepsilon }\label{yahdi1}
.\end{equation}\newline $\bullet$ if $\rho^{0}$ is of short range
type, we have
\begin{equation} \liminf_{\varepsilon \to 0^{+}}
\frac{\log\Big|\log \Big ( N(E_+
 +\varepsilon)-N(E_+)\Big)\Big|}{\log\varepsilon}\geq \liminf_{\varepsilon \to 0^{+}}
\frac{\log\Big|\log  (
\mathbb{P}_{\varepsilon,\alpha,2})\Big|}{\log\varepsilon}\label{yahdi2}
.\end{equation} Now one deals with the estimation of
$\mathbb{P}_{\varepsilon,\alpha,1}$ and
$\mathbb{P}_{\varepsilon,\alpha,2}$. We start by:\newline
$\bullet$ The estimation of
$\mathbb{P}_{\varepsilon,\alpha,1}$.\newline Let $\nu \in
(d,d+2]$. First we notice that if
$$ \omega_{\gamma}\leq \varepsilon^{1+\alpha}\ \ \text{for}\ \ |\gamma|\leq \varepsilon^{-(1-\alpha)/2 }$$
and $$\omega_{\gamma}\leq
\varepsilon^{1+\alpha}\Big(1+dist(\gamma,C_{0,\varepsilon^{-(1-\alpha)/2}})\Big)^{(\nu-d)(1-\alpha)}\
\ \text{for}\ \
 \varepsilon^{-(1-\alpha)/2}<|\gamma|\leq \varepsilon^{-\frac{1+2\alpha}{\nu-d}},$$ then
$$ \sum_{\gamma\in \mathbb{Z}^d}\omega_{\gamma}(1+|\beta-\gamma|)^{-\nu}\leq \varepsilon^{1+\alpha}.$$
So \begin{equation} \mathbb{P}_{\varepsilon,\alpha,1}\geq
\mathbb{P}_{2}\cdot \mathbb{P}_{1}. \label{prod1}\end{equation}
Where
\begin{equation*}
\mathbb{P}_{1}=\mathbb{P}\{\omega;\ \forall\  \gamma\  \text{such
that }\ |\gamma|\leq \varepsilon^{-(1-\alpha)/2},\
\omega_{\gamma}\leq \varepsilon^{1+\alpha}\},
\end{equation*} and
\begin{multline*}\mathbb{P}_{2}=\mathbb{P}\Big\{\omega;\ \forall\ \
\gamma\  \text{such that } \
\varepsilon^{-(1-\alpha)/2}<|\gamma|\leq
\varepsilon^{-\frac{1+2\alpha}{\nu-d}},\\ \omega_{\gamma}\leq
\varepsilon^{1+\alpha}\Big(1+dist(\gamma,C_{0,\varepsilon^{-(1-\alpha)/2}})\Big)^{(\nu-d)(1-\alpha)}\Big\}.
\end{multline*} As the random variables are i.i.d we
get that
\begin{equation}
\mathbb{P}_{1}=\Big(\mathbb{P}\{\omega_{0}\leq
\varepsilon^{1+\alpha}\}\Big)^{\varepsilon^{-d(1-\alpha)/2}}\label{prod2}
\end{equation}
and
\begin{equation}\mathbb{P}_{2}=\prod_{\varepsilon^{(1-\alpha)/2}<|\gamma|\leq
\varepsilon^{-\frac{1+2\alpha}{\nu-d}}}\mathbb{P}\Big(\omega_0\leq
 \varepsilon^{1+\alpha}(1+dist(\gamma,C_{0,\varepsilon^{-(1-\alpha)/2}}))^{(\nu-d)(1-\alpha)})\Big).
\label{prod3}\end{equation} Now by applaying the logarithm to
(\ref{prod1}) and taking into akount (\ref{prod2}) and
(\ref{prod3}) wile using (\ref{you1}), for $\alpha $ and $
\varepsilon$ small enough we get that\newline
\begin{multline}
\log \mathbb{P}_{\varepsilon,\alpha,1} \geq\\
-\varepsilon^{-(\kappa+d/2)(1+\alpha)}-\varepsilon^{-\kappa(1+\alpha)}\sum_{\varepsilon^{-(1-\alpha)/2}\leq
|\gamma|\leq \varepsilon^{-\frac{1+2\alpha}{\nu-d}}}
\Big(1+dist(\gamma,C_{0,\varepsilon^{-(1-\alpha)/2}})\Big)^{-\kappa(\nu-d)(1-\alpha)}\label{you3}
.\end{multline} As if $(\nu-d)\kappa>d$ the sum in (\ref{you3})
converges when $\alpha$ is chosen small enough such that
$(1-\alpha)(\nu-d)\kappa>d$. So we get,
\begin{equation}
\liminf_{\varepsilon \to 0^{+}} \frac{\log\Big|\log  (
\mathbb{P}_{\varepsilon,\alpha,2})\Big|}{\log\varepsilon}\geq
-(1+\alpha)\Big(\kappa +\frac{d}{2}\Big).
\label{mouch5}\end{equation} In the case when $(\nu-d)\kappa<d$,
for $\varepsilon$ small one computes the sum in (\ref{you3}) we
get the following estimation
\begin{multline}
\sum_{\varepsilon^{-(1-\alpha)/2}\leq |\gamma|\leq
\varepsilon^{-\frac{(1+2\alpha)}{(\nu-d)}}}
\Big(1+dist(\gamma,C_{0},\varepsilon^{-(1-\alpha)/2})\Big)^{-\kappa(\nu-d)(1-\alpha)}\leq
\\ C\cdot  \varepsilon ^{\kappa(1-\alpha)}\cdot
\varepsilon^{-d(1+\alpha)/(\nu-d)}. \label{you4}\end{multline}
Using equations (\ref{you3})  and (\ref{you4}) and the fact that $
\frac{d}{\nu-d}-\kappa+\frac{d}{2}\geq 0$ we get
$$
\log \mathbb{P}_{\varepsilon, \alpha,1}\geq -C\varepsilon ^{-(\kappa+d/2)(1+\alpha)}\cdot \varepsilon^{-(d(1+\alpha)-2
\alpha\kappa)/(\nu -d)}.
$$
We apply the logarithm  into the last equation taking into a count
(\ref{you1}), (\ref{yahdi1}) and (\ref{mouch5}) we get
(\ref{hifr}).\newline $\bullet$ The estimation of
$\mathbb{P}_{\varepsilon,\alpha,2}$. Let us notice that there
exists $C>0$ such that we have
$$\mathbb{P}_{\varepsilon,\alpha,2}\geq \mathbb{P}\Big \{\omega; \forall \gamma \in \Lambda_{\alpha}(\varepsilon^{s});\omega_{\gamma}\leq \frac{\varepsilon^{1+\alpha}}{C} \Big \}.$$
As the random variables are  i.i.d we get that
$$
\mathbb{P}_{\varepsilon,\alpha,2}\geq
\Pi_{\gamma\in\Lambda_{\alpha}(\varepsilon^{s})}\mathbb{P}\{\omega;
\omega_{\gamma} \leq
\frac{\varepsilon^{1+\alpha}}{C}\}=\Big(\mathbb{P}\{\omega_{0}\leq
\frac{\varepsilon^{1+\alpha}}{C}\}\Big)^{\sharp
\Lambda_{\alpha}(\varepsilon^{s})}.
$$
Now taking into account (\ref{you1}), (\ref{yahdi2}), the fact
that $\sharp
\Lambda_{\alpha}(\varepsilon)=\varepsilon^{-d(\frac{1}{2}+\alpha)}$,
and $\alpha>0$ small we end the proof of (\ref{mouch1}). So the
proof of
Theorem \ref{LL1} is ended. \hfill $\Box$
\subsection{The upper bound}
To prove the upper bound, we compare $N(E_{+}+ \varepsilon)-N(E)$
to the IDS some reduced operators. More precisely,
we prove that for an energy $E$ close to $E_+$, $%
N(E)-N(E_+)$ can be upper bounded by the IDS of some random and
bounded operator. Indeed we have,
\begin{lem}\text{\cite{jmp}}
\label{bb22}Let $H_{\omega}$ be the operator defined by (\ref{b}). We assume that \textbf{%
(A.1)},{\ \textbf{(A.2)}}  and {\textbf{(A.3)}} hold. There exists
$E_{0}>E_+$ and $C>1$ such that, for $E_+\leq E\leq E_{0}$ we have
\begin{equation}
0\leq N(E)-N(E_+)\leq N_{{{\mathcal{E}}_{0}}}\Big(C\cdot
(E-E_+)+E_+\Big) \label{T}
\end{equation}
where $N_{{\mathcal{E}}_{0}}$ is the IDS of $H_{\omega}^0=\Pi_{0}H_{\omega}%
\Pi_{0}$ and $\Pi_{0}$, is the spectral projection for $H_0 $ on
the band starting at $E_{+}$.
\end{lem}
{\bf{Proof}}: The proof of Lemma \ref{bb22} is given in the case
of acoustic operators in \cite{jmp}. It is still true in the
divergence case. It is based on a localization in energy for the
density of states. It goes as follows: We approach the density of
states of $H_{\omega}$ by the density of states of periodic
approximations, see section \ref{abour}. In a neighborhood of
$E_+$, we control the behavior of the density of states of
periodic approximations via the density of states of periodic
approximations of the reference operator i.e
$H_{\omega}^{0}=H_0^0+V_{\omega}^0$. We then compute the limit for
the density of states of the reference operators and we obtain the
sought for result.$\hfill \Box $
\subsubsection{The short range case} \label{6.1} We recall that in this case we
assume that the IDS, $n$ of the background operator $H_0$ is
non-degenerate. We prove the following theorem:
\begin{theo}
\label{LL12}Let $H_{\omega}$ be the operator defined by (\ref{b}).
We assume that \textbf{(A.1)} and \textbf{(A.2)} hold  and $n$ is non-degenerate at $E_{+}$,
then
$$
\limsup_{\varepsilon \longrightarrow 0^{+}}\frac{\log | \log
(N(E_{+}+\varepsilon )-N(E_{+}))|}{\log \varepsilon} \leq
-(\frac{d}{2}+\kappa).
$$
\end{theo}
\textbf{Proof:} Let us notice that by Lemma \ref{bb22} to prove
Theorem \ref{LL12} it suffices to get the same upper bound for the
reference operator. This represents several advantages: first, $%
H_{\omega}^{0}$ it is a bounded random operator and
 equivalent to a random Jacobi matrix acting on $L^{d}(\Bbb{T}^*)\otimes {\Bbb{C}}^{n_0}$
 (Here $n_0$ is the number of Floquet eigenvalues generating the band starting in $E_+$).
 The second advantage is that
while, $E_+$ is an interior edge of a gap for $H_{\omega}$, it
becomes the bottom of the spectrum for $H_{\omega}^{0}$.\newline
The idea of the proof is based on the uncertainly principle.
Indeed, as $V_{\omega}^0\geq 0$, if a vector minimizes
$H_{\omega}^0$, it necessarily minimizes $H^0 _ {0}=\Pi _ {0}
H_{0} \Pi _ {0} $; hence, it has to be concentrated in the
quasimomentum $\theta$ near the zeros of
$(E_{j}(\theta)-E_{+})_{1\leq j\leq n_{0}}$. For this, we have to
take into account all the points where the Floquet eigenvalues
reach $E_+$. Let $\theta^0$ be one of those points. As $n$ is
non-degenerate, for $c>0$ small, in a neighborhood of $\theta^0$
we have \cite{jmp}
\begin{equation}
\mathbb{D}=c\sum_{j=1}^{d}(1-\cos (\theta_j-\theta_j^0))\leq
H_0^0-E_+\cdot I_{d}.
\end{equation}
Here $\mathbb{D}$ is acting on $L^{2}(\mathbb{T}^*)\otimes
\mathbb{C}^{n_{0}}$ and $I_{d}$ is the identity matrix.\newline
We recall that $V_{\omega}$ is the operator defined by
\begin{equation}
V_{\omega}=-\nabla\Big( \sum_{\gamma\in
\mathbb{Z}^d}\omega_{\gamma}\rho^{0}(\cdot -\gamma)\Big)\nabla.
\end{equation}
It is proved \cite{jmp} that $V^{0}_{\omega} $ can be lower bounded
by $$\displaystyle V_{2,\omega} ^a=\sum_{\gamma \in
{\Bbb{Z}^d}}\omega _{\gamma}\Pi _{\gamma}.$$ Here $ \Pi _{\gamma}$
is the orthogonal projection on the vector $ \theta\mapsto
e^{i{\gamma} \theta}$ in
$L^2({\mathbb{T}^*})\otimes\mathbb{C}^{n_{0}}$. Now using the following
(unitary operator) discrete Fourier transformation defined from $l^{2}(\mathbb{Z}^d)$ to $L^{2}([0,2\pi]^{d})$ by
$$
\mathcal{F}(u)(k)=\widehat{u}(k)=\sum_{n\in
\mathbb{Z^d}}u(n)e^{-in\cdot. k},
$$
we get that $\mathbb{D}$ is unitarly equivalent to the usual
discrete Schr\"odinger operator. So $H_{\omega}^0$ is lower
bounded by  some opertaor which it self unitarly equivalent to the usual discrete random operator whose behavior of
the IDS at the edges of the spectral gaps is already known
\cite{kp1,jmp}. This lower bound on the operator immediately
yields an upper bound on the density of states.
\subsubsection{The long range case }In this section we shale prove: \begin{theo}\label{adel}
Let $H_{\omega}$ be the operator defined by (\ref{b}). We assume
that {\bf{(A.1)}}, {\bf{(A.2)}} hold.\\ If $\frac{d}{\nu-d}>
\kappa +\frac{d}{2}$ then,
\begin{equation}\limsup_{\varepsilon \to 0^+} \frac{\log|\log \big( N(E_+
+\varepsilon)-N(E_+)\big)|}{\log\varepsilon}\leq
-\frac{d}{\nu-d}.\label{mouch3}
\end{equation}
If $n$, the IDS of $H_0$ is non-degenerate then,
\begin{equation}
\limsup_{\varepsilon \to 0^+} \frac{\log|\log \big( N(E_+
+\varepsilon)-N(E_+)\big)|}{\log\varepsilon}\leq
-\sup(\frac{d}{2}+\kappa,\frac{d}{\nu-d}).\label{mouch2}
\end{equation}
\end{theo}
{\bf{Proof:}}\newline $\bullet$ If $\frac{d}{\nu-d}> \kappa
+\frac{d}{2}$.\newline Notice that in this case we have no
assumption made on the behavior of $n$, the IDS of the periodic
operator. The proof goes exactly as the one given in \cite{long},
for this we omit details. From Lemma \ref{mez} and for
$\eta_{0}>1/(\nu-d)$ and $k\sim \varepsilon ^{-\delta}$ such that
$\delta>\nu_{0}$ the proof of (\ref{mouch3}) is reduced to prove
that
\begin{equation}\label{bel}
\limsup_{\varepsilon\rightarrow
0^{+}}\frac{\log\Big|\log\Big(\mathbb{E}(N_{\omega,k}
(E_{+}+\varepsilon)-N_{\omega,k}(E_{+}))\Big)\Big|}{\log\varepsilon}\leq
-\frac{d}{\nu-d}.
\end{equation}
\begin{lem}\cite{jmp}\label{3.1}Let $k\sim\varepsilon^{-\rho}$ with
$\rho>1/(\nu-d)$. Define the event,
$${\bf{E}}_{\varepsilon,\omega}=\Big\{\omega; V_{\omega,k}\geq -\varepsilon
\Delta=-\varepsilon\sum_{i=1}^{d}\partial_{x_{i}}^{2}\Big\}.$$
Here we recall that $$ V_{\omega,k}=\nabla \cdot \Big(\sum_{\gamma
\in C_{k}\cap {\mathbb{Z}^{d}} }\omega _{\gamma }\sum_{\beta \in
(2k+1){\mathbb{Z}^{d}}}\rho^{0}(\cdot -\gamma -\beta)\Big)\cdot
\nabla.$$Then ${\bf{E}}_{\varepsilon, \omega}$ has a probability
at least $1-\mathbb{P}_{\varepsilon}$ where
$\mathbb{P}_{\varepsilon}$ satisfies
\begin{equation}\label{3.2}
\limsup_{\varepsilon\rightarrow
0^{+}}\frac{\log|\log(\mathbb{P}_{\varepsilon})|}{\log
\varepsilon}\leq -\frac{d}{\nu-d}.
\end{equation}
\end{lem}
Using the fact that if for some $C>0$ (depending only on $\delta$
and $\rho^*$), $V_{\omega,k}\geq -C\varepsilon \Delta $, then the
spectrum of $H_{\omega,k}$ does not intersect
$(E_{+},E_{+}+\varepsilon)$ for $\varepsilon$ small.\newline One
computes
\begin{eqnarray*}
\mathbb{E}\Big(N_{\omega,k}(E_{+}+\varepsilon)-N_{\omega,k}(E_{+})\Big)&=&
\mathbb{E}\Big([N_{\omega,k}(E_{+}+\varepsilon)-N_{\omega,k}(E_{+})]_{\mathbf{1}_{\{\omega;V_{\omega,k}\geq
-C\varepsilon \Delta\}}}\Big)\\
&+&\mathbb{E}\Big([N_{\omega,k}(E_{+}+\varepsilon)-N_{\omega,k}(E_{+})]_{\mathbf{1}_{\{\omega;V
_{\omega,k}<
-C\varepsilon \Delta\}}}\Big) \\
&\leq&C\mathbb{P}(\{\omega;V_{\omega,k}<-C\varepsilon
\Delta\})\\
&=& C(1-\mathbb{P}({\bf{E}}_{C\cdot
\varepsilon,\omega}))=C\mathbb{P}_{C\cdot \varepsilon}.
\end{eqnarray*}
Here, we have used the fact that $N_{\omega,k}$ is bounded,
locally uniformly in energy, uniformly in $\omega$, $k$ by $C$.
Taking (\ref{3.2}) into account, we end the proof of (\ref{bel})
and so (\ref{mouch3}) is proved. $\Box$\newline Now one deals with
the proof of (\ref{mouch2}). We recall that here one supposes once
more that $n$ is non-degenerate. The idea is similar to the short
range case (we will compare the IDS of our operator to the IDS of
another one) and
 we will follow and use results given in \cite{Pasfig, sim}.
 Let $N^{a}$ be the IDS of the follwoing Anderson discrete
operator acting on $l^{2}(\mathbb{Z}^d)$:
\begin{equation}
 (H_{\omega}^au)(\alpha)= E_+\cdot
 u(\alpha)+\sum_{|\alpha-\beta|=1}(u(\alpha)-u(\beta))+(V_{\omega}^{a}u)(\alpha).\label{mouch6}
 \end{equation}
 Here $V^a_{\omega}$  the diagonal infinite matrix with $v_{\alpha}(\omega)=\sum_{\beta\in\mathbb{Z}^d}\omega_{\beta}
 (1+|\alpha-\beta|)^{-\nu}$ for the $\alpha^{\text{th}}$ diagonal coefficient.\\
For $k\in \mathbb{N}^*$ and $u\in l^{2}(\mathbb{Z}^d\cap C _{k})$,
let $H_{0}^{k},V_{\omega}^{k}$ and $H_{\omega}^{k}$ be the
following discrete operators
$$
(H_{0}^{k}u)(\alpha)=E_+\cdot u(\alpha)+\sum_{|\alpha-\beta|=1,
\beta\in C_{k}}(u(\alpha)-u(\beta)), \ \
(V_{\omega}^{k}u)(\alpha)=v_{\alpha}(\omega)u(\alpha)$$ and
$$
 H_{\omega}^{k}=H_{0}^{k}+V_{\omega}^{k}.
$$ Let $N_{k}^{a}$ be the IDS of $H_{\omega}^{k}$ defined by:
$$
N_{k}^{a}(E)=\frac{1}{(2k+1)^d}\cdot \mathbb{E}\Big(\sharp\{
\text{eigenvalues of }\ H_{\omega}^{k}\ \text{less or equal to }\
E\} \Big).
$$
From \cite{Pasfig,sim}, we know that for a good choice of $k$, the
IDS at
energy $E$ is quite well approximated by the probability to
find a state energy less than $E$.  Precisely we have the
following relation:
\begin{equation}
N^{a}(E)\leq N_{k}^{a}(E)\leq C\cdot \mathbb{P}_{k}(E).
\end{equation}
Here
$$
\mathbb{P}_{k}(E)=\mathbb{P}\Big(\Big\{H_{\omega}^{k}\ \text{admits at least an eigenvalue less than}\ E\Big\}\Big).
$$
To estimate this probability  one proceeds as previously and lower bound $H_{\omega}^{k}$ by
 $H_{\widetilde{\omega}}^{k}$; obtained
for $\delta>0$ by changing $\omega_{\gamma}$ by
$\widetilde{\omega}_{\gamma}=\omega_{\gamma}$ if
$\omega_{\gamma}\leq \delta$ and
$\widetilde{\omega}_{\gamma}=\delta$ if not. So if we denote the
IDS of $H_{\widetilde{\omega}}^{k}$ by $\widetilde{N}_{k}^{a}$
then we have,
$$
N_{k}^{a}(E)\leq \widetilde{N}_{k}^{a}(E).
$$
One takes $k=c(E-E_{+})^{-\frac{1}{2}}$ and $\delta=(E-E_+)/c$ positives. For a good choice of $c$ and thus of $\delta$, and $E$ in a neighborhood
 of $E_+$, we get
\begin{equation}
\mathbb{P}_{k}(E)\leq \mathbb{P}\Big(\Big\{\widetilde{\omega};\ \frac{1}{(2k+1)^d}\sum_{|\gamma|\leq k}v_{\gamma}(\widetilde{\omega})\leq \delta/K\Big\}\Big).
\end{equation}
Now we have to estimate the last probability. Let $0<\alpha <1$, for some $C>1$ we have:
$$
\frac{1}{(2k+1)^d}\sum_{|\gamma|\leq
k}v_{\gamma}(\widetilde{\omega})= \frac{1}{(2k+1)^d}\sum_{\beta\in
\mathbb{Z}^d}\widetilde{\omega}_{\beta}\Big(\sum_{|\gamma|\leq
k}(1+|\gamma-\beta|)^{-\nu}\Big).
$$
$$
\geq \frac{1}{C(2k+1)^d}\sum_{|\beta|\leq
k}\widetilde{\omega}_{\beta}+\frac{1}{C} \sum_{k<|\beta|\leq
\delta^{-(\nu-d)(1-\alpha)}}\widetilde{\omega}_{\beta}(1+|\beta|+k)^{-\nu}.$$
Thus
\begin{multline}
\mathbb{P}\Big(\Big\{\widetilde{\omega};\ \frac{1}{(2k+1)^d}\sum_{|\gamma|\leq k}
v_{\gamma}(\widetilde{\omega})\leq \delta/K\Big\}\Big)\leq \\
\mathbb{P}\Big(\Big\{\widetilde{\omega};\
\frac{1}{C(2k+1)^d}\sum_{|\beta|\leq
k}\widetilde{\omega}_{\beta}\leq \delta/K\ \text{and}\ \
\frac{1}{C} \sum_{k<|\beta|\leq
\delta^{-(\nu-d)(1-\alpha)}}\widetilde{\omega}_{\beta}(1+|\beta|+k)^{-\nu}\leq
\delta/K\Big\}\Big).
\end{multline}
As the random variables are i.i.d we get for
$$
\mathbb{P}_{1}=\mathbb{P}\Big(\Big\{\widetilde{\omega};\ \frac{1}{C(2k+1)^d}
\sum_{|\beta|\leq k}\widetilde{\omega}_{\beta}
\leq \delta/K\ \Big\}\Big)
$$
and
$$
\mathbb{P}_{2}=\mathbb{P}\Big(\Big\{\widetilde{\omega};
\frac{1}{C} \sum_{k<|\beta|\leq
\delta^{-(\nu-d)(1-\alpha)}}\widetilde{\omega}_{\beta}(1+|\beta|+k)^{-\nu}\leq
\delta/K\Big\}\Big),
$$
we have
$$
\mathbb{P}\Big(\Big\{\widetilde{\omega};\ \frac{1}{(2k+1)^d}\sum_{|\gamma|\leq k}v_{\gamma}(\widetilde{\omega})\leq \delta/K\Big\}\Big)\leq \mathbb{P}_{1}\cdot \mathbb{P}_{2}.
$$
The estimation of $\mathbb{P}_{1}$ and $\mathbb{P}_{2}$ is based
on large deviation results \cite{da}. Briefly the idea is the
following. Let $t>0$. Using the Markov inequality one estimates
\begin{equation}
\mathbb{P}_{1} \leq
\mathbb{E}\Big(e^{t(\delta/K-\frac{1}{(2k+1)^d}\sum_{|\beta|\leq
k}\widetilde{\omega}_{\gamma})}\Big). \label{youa1}\end{equation}
As the random variables are i.i.d one gets
\begin{eqnarray}
\mathbb{E}\Big(e^{t(\delta/K-\frac{1}{(2k+1)^d}\sum_{|\beta|\leq
k}\widetilde{\omega}_{\gamma}})\Big)&=&e^{t\delta/K}\prod_{|\beta|\leq
k}
\mathbb{E}\Big(e^{-t\widetilde{\omega}_{\beta}/(2k+1)^{d}}\Big)\\
&=&e^{t\delta/K}e^{(2k)^d\log\mathbb{E}(e^{-t\widetilde{\omega}_{0}/(2k+1)^d})}\label{youwa2}
\end{eqnarray}
For $k$ big enough using Taylor expansion  of $e^{t}$, we get that
$$
\mathbb{E}\Big(e^{-t\frac{\widetilde{\omega
}_{0}}{(2k)^d}}\Big)=1-\frac{t\mathbb{E}(\widetilde{\omega}_{0})}{(2k+1)^d}+
o(\frac{t^{2}}{(2k+1)^{2d}}).
$$
So we get
$$
\mathbb{P}_1\leq
e^{t\delta/K}e^{(2k)^d\log\Big(1-\frac{t\mathbb{E}(\widetilde{\omega}_{0}=0)}{(2k+1)^d}+
o(\frac{t^{2}}{(2k+1)^{2d}}) \Big)}.
$$
As $\mathbb{E}(\widetilde{\omega}_{0})>0$ for $k$ big enough and
$t>0$ well chosen we get that there exists $C>0$ such that
 $$
\mathbb{P}_{1}\leq e^{\Big((2k)^d\log
(\mathbb{P}\{\widetilde{\omega}_{0}=0 \}\Big)/C}.$$ The same
computation as above gives
$$
\mathbb{P}_{2}\leq e^{-\Big(k^d|\log
(\mathbb{P}\{\widetilde{\omega}_{0}\leq
\delta\})|+\delta^{-d/(\nu-d)}\Big)/C}.
$$
So, we get that:
$$
\log N^{a}(E)\leq -\Big(k^d|\log
\mathbb{P}(\{\widetilde{\omega}_{0}=0\})|+|\log
(\mathbb{P}(\widetilde{\omega}_{0}=0))|+\delta^{-d/(\nu-d)}\Big)/C.
$$
Now, for  $k=c(E-E_+)^{-\frac{1}{2}}$ and $\delta=(E-E_{+})/c$ we
get (\ref{mouch2}) for $N^{a}$.\newline  As $n$ is non-degenerate
we get that (see the short range case and \cite{jmp})
$H_{\omega}^a$ is unitarly equivalent by means of the Fourier
transformation to some operator which lower bound $H_{\omega}^0$.
This yelds (\ref{mouch2}) for $N$.\hfill $\Box$.
\newline
\textit{$\mathbf{Acknowledgements.}$ I would like to thank\
professor Mabrouk Ben Ammar and my colleague Adel Khalfallah for
theirs supports.}

\end{document}